\documentclass{amsart}
\usepackage[utf8]{inputenc}
\usepackage{amsmath, amsthm, amssymb, amsfonts}
\usepackage{bm}
\usepackage{hyperref}
\usepackage{dsfont, comment}
\usepackage[utf8]{inputenc}
\usepackage{rotate}

\usepackage{tikz}
\usepackage{tikz-cd}
\usepackage[arrow,matrix,curve]{xy} 	
\usepackage{xcolor}
\usepackage{todonotes}
\usepackage{alltt} 
\usepackage{amsmath,amssymb}
\usepackage{url}
\usepackage{todonotes}
\usepackage[margin=1in]{geometry}

\theoremstyle{plain}
\newtheorem{theorem}{Theorem}[section]

\theoremstyle{definition}

\newtheorem{remark}[theorem]{Remark}
\newtheorem{example}[theorem]{Example}

\newcommand{\norm}[1]{\left\lVert#1\right\rVert}

\definecolor{darkpastelgreen}{rgb}{0.01, 0.75, 0.24}

\title{Sharp Bounds for Limits of $N$-Point Correlation}
\author{Christian Wei\ss{}}
\date{\today}

\begin{document}

\begin{abstract} We consider the question which functions $f$ may appear as limiting functions of the $N$-point correlation. We obtain sharp point-wise lower and upper bounds. The explicit examples constructed to show the sharpness of the bounds are based on van der Corput sequences.
\end{abstract}

\maketitle

\section{Introduction}

For $N \in \mathbb{N}$ and $s \geq 0$, the $N$-point correlation of a sequence $(x_n)_{n \in \mathbb{N}} \subset [0,1)$ is defined by
\[
R_{N}(s) = \frac{1}{N} \# \left\{ 1 \leq i \neq j \leq N \, : \, \norm{x_i-x_j} \leq \frac{s}{N} \right\},
\]
where $\norm{\cdot}$ denotes the distance to the nearest integer. It is well-known that  for a sequence of random variables $(X_n)_{n \in \mathbb{N}}$ drawn independently from uniform distribution $R_N(s)$ almost surely converges to $2s$ for all $s \geq 0$. If this limit relation holds, the sequence $(x_n)_{n \in \mathbb{N}}$ has Poissonian pair correlation, compare \cite{ALP18}. Although this is a generic property of independent, uniformly distributed random variables, it is a difficult task to find \textit{deterministic} sequences possessing Poissonian pair correlation. Only relatively few classes of examples have been found, see e.g. \cite{BMV15,LST24,LT22}.\\[12pt]
On the contrary, the case where $\lim_{N \to \infty} R_{N}(s)= f(s)$ exists for all $s \geq 0$ but $f(s) \neq 2s$, has been much less analyzed in the literature. We say that the sequence $(x_n)_{n \in \mathbb{N}}$ has $f$-pair correlation in this case, compare \cite{FW25}. The present paper approaches the question which functions $f$ may appear as limiting functions. While it is immediate that $f(s)$ must be a non-decreasing function with $f(s) \to \infty$ and that  $f(0) = 0$ if the number of multiply appearing elements in the sequence is of order $o\left(N^{1/2}\right)$, only little else is known: there are some examples in the non-generic case, where the limiting functions have been explicitly calculated, see  \cite{ALP18, Lut20, MS13, Say23, Wei23}. Moreover, it has been proven in \cite{FW25} that $f(s)$ cannot be continuous if the number of different gap lengths, i.e. the distances between geometrically neighboring elements of the finite sequence $x_1,\ldots,x_N$, is uniformly bounded. To the best of the author's knowledge no other general property of $f(s)$ has been shown in the literature. Here, we prove the following sharp bounds which hold for all possible limit functions $f(s)$.

\begin{theorem} \label{thm:main}
Let $(x_n)_{n \in \mathbb{N}}\subset \mathbb [0,1)$.
\begin{enumerate}
\item[(i)] For every $s\geq 0$, we have
\[
\liminf_{N \to \infty} R_N(s) \geq G(s),
\]
where $G(0)=0$, and for $s>0$, we define
\[
G(s) = \frac{h(2s-h-1)}{s}, \qquad h=\lfloor s\rfloor .
\]
Equivalently, if $s=h+\theta$ with $h=\lfloor s\rfloor$ and $\theta\in[0,1)$, then
\[
G(s) = s-1+\frac{\theta(1-\theta)}{s}.
\]
\item[(ii)] If the limit function $f: \mathbb{R}^+_0 \to \mathbb{R}^+_0$ exists, then for every $0<a<b$, it holds that
\[
f(b)\leq C(a,b)\left(1+f(a)\right),
\]
where
\[
C(a,b) = 2\left( \left\lfloor \frac{b}{a} \right\rfloor +2 \right) + 1.
\]
\item[(iii)] If $a>0$ and
\[
f(a)=0,
\]
then for every $b>0$,
\[
f(b)\leq 2\left\lfloor \frac{b}{a}\right\rfloor.
\]
Equivalently, if
\[
\lim_{N \to \infty} R_N(a)= 0,
\]
then
\[
\limsup_{N\to\infty} R_N(b) \leq 2\left\lfloor \frac{b}{a}\right\rfloor .
\]
\item[(iv)] The bound in (i) is sharp while the bound in (iii) is sharp in the sense that for every integer $k \geq 1$ and every $a \in (k/(k+1),1)$, there exists a sequence $(x_n)_{n \in \mathbb{N}}$ such that $\liminf_{N \to \infty} R_N(a) = 0$, but $\limsup_{N \to \infty} R_N(h) \geq 2k = 2 \left\lfloor \frac{k}{a} \right\rfloor$.
\end{enumerate}
\end{theorem}
In particular, it follows that if the limit function $f(s)$ exists for all $s \geq 0$ and $f(a) = 0$ for some $a>0$, then 
\[
\max(0,s-1) \leq f(s) \leq 2 \left\lfloor \frac{s}{a} \right\rfloor,
\]
for all $s \geq a$. If $f(s) \neq 0$ for all $s > 0$, then we only obtain the weaker bounds
\[
\max(0,s-1) \leq f(s) \leq (1+f(1)) (2 \left\lfloor s \right\rfloor+5).
\]
In other words, our results provide universal necessary conditions, including sharp point-wise lower and upper bounds for a function to arise as a limit of the $N$-point correlation function. A comparable result for the weak $N$-point correlation, which excluded the case addressed here, was derived recently in \cite{Wei26d}.\\[12pt]
The proof of the lower bound relies on dividing $[0,1)$ into $M_N = \lceil \frac{N}{s} \rceil$ intervals of equal length and some point counting arguments. Similarly, the upper bounds follow by considering $L_N = \lceil \frac{N}{a} \rceil$ intervals of equal length. A comparable approach was applied in several proofs of \cite{HZ24} to analyze weak pair correlation. The explicit examples constructed to show sharpness of the bounds are based on van der Corput sequences. While a van der Corput sequence suffices to show the sharpness of (iii), to prove the sharpness of (i) we shift it by a small (in our case deterministic) error, an idea also applied in \cite{Por26}. 
\begin{example} In \cite{Wei23}, the limit function $f_1(s)$ of the (low-dispersion) sequence $x_n := \left\{ \frac{\log(2n-1)}{\log(2)} \right\}$, where $\{ \cdot \}$ denotes the fractional part of a number, was calculated as
\[ 
f_1(s) = \begin{cases}  
0 & 0 \leq s < \log(4)^{-1}\\
3k - \frac{3k^2+2k}{4\log(2)s} & k \log(4)^{-1} \leq s < (k+1) \log(4)^{-1}, k \in \mathbb{N} \ \textrm{even},\\ 3k + 1 - \frac{3k^2+4k+1}{4\log(2)s} & k \log(4)^{-1} \leq s < (k+1) \log(4)^{-1}, k \in \mathbb{N} \ \textrm{odd}.\end{cases}
\]
In Figure~\ref{thm:main}, we compare this function to the bounds calculated in Theorem~\ref{thm:main} (i) and (iii) when choosing $a=\log(4)^{-1}$. We see that the function $f_1(s)$ is indeed within the range defined by the bounds.
\end{example}
\begin{figure}[ht]
    \centering
    \includegraphics[width=0.5\linewidth]{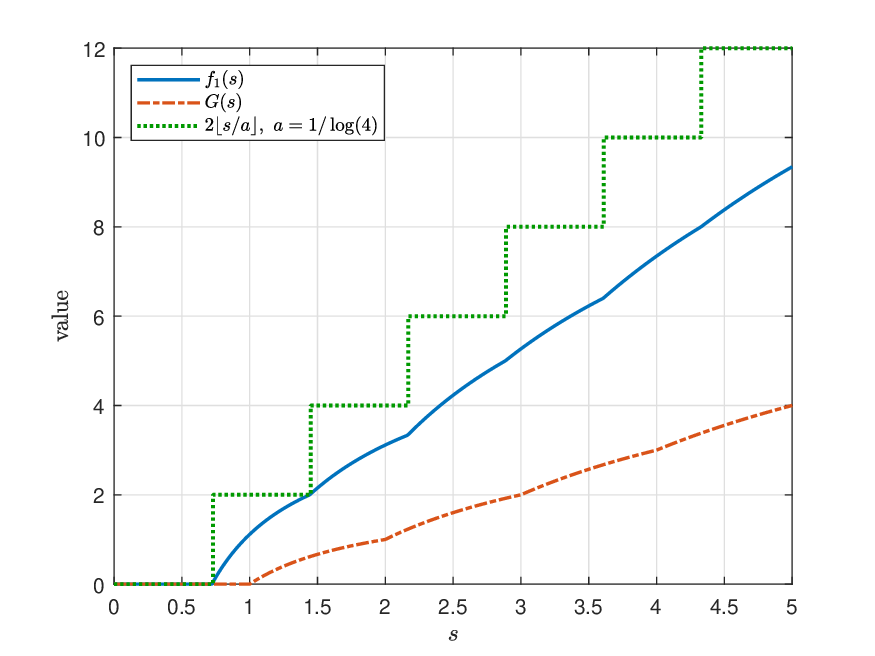}
    \caption{Comparison of the limit function for the sequence $x_n := \left\{ \frac{\log(2n-1)}{\log(2)} \right\}$ with the bounds from Theorem~\ref{thm:main}.}
    \label{fig:placeholder}
\end{figure}
\section{Proof of the results}
For $s\geq 0$, set
\[
A_N(s) = \#\left\{ 1\leq i\neq j\leq N\, : \, \norm{x_i-x_j}\leq \frac{s}{N} \right\}.
\]
(i) For $s=0$, the claim is simply $\liminf R_N(s)\geq 0$, so there is nothing to prove. Now fix $s>0$ and partition $[0,1)$ into 
\[
M_N:=\left\lceil \frac{N}{s}\right\rceil
\]
intervals $I_1,\dots,I_{M_N}$ of equal length $1/M_N$. Since $\frac{1}{M_N}\leq \frac{s}{N}$, any two points lying in the same interval are at distance at most $s/N$. Setting
\[
n_m=\#\{1\leq i\leq N:x_i\in I_m\}.
\]
we observe
\[
A_N(s)\geq \sum_{m=1}^{M_N}n_m(n_m-1).
\]
Next, we minimize the right-hand side under the constraint $\sum_{m=1}^{M_n} n_m=N$. Uniquely writing
\[
N=q_N M_N+r_N,
\]
for some $0\leq r_N<M_N$, we see that the minimum is attained when $r_N$ of the $n_m$ are equal to $q_N+1$, and the remaining $M_N-r_N$ ones are equal to $q_N$. It follows that
\[
\sum_{m=1}^{M_N}n_m(n_m-1) \geq (M_N-r_N)q_N(q_N-1)+r_N(q_N+1)q_N.
\]
Therefore,
\[
R_N(s) \geq \frac{M_N q_N(q_N-1)+2q_N r_N}{N}.
\]
Since $t_N:=\frac{N}{M_N}$ converges to $s$ as $N \to \infty$ and $q_N=\left\lfloor t_N\right\rfloor$ as well as $\frac{r_N}{M_N}=t_N-q_N$, we may rewrite the right-hand side as
\[
\frac{q_N(q_N-1)+2q_N(t_N-q_N)}{t_N}.
\]
or equivalently,
\[
R_N(s)
\geq
\frac{q_N(2t_N-q_N-1)}{t_N}.
\]
If $s\notin\mathbb N$, then $q_N=\lfloor t_N\rfloor\to \lfloor s\rfloor=h$, and therefore
\begin{align} \label{lb:RNs}
\liminf_{N\to\infty} R_N(s) \geq \frac{h(2s-h-1)}{s}.
\end{align}
If $s\in\mathbb N$, then $t_N\to s$, and for all sufficiently large $N$
we have
\[
q_N=\lfloor t_N\rfloor\in\{s-1,s\}.
\]
In both cases the right-hand side of \eqref{lb:RNs} converges to $s-1$ as $N \to \infty$.
Thus, we obtain
\[
\liminf_{N\to\infty}R_N(s) \geq \frac{h(2s-h-1)}{s}, \qquad h=\lfloor s\rfloor.
\]
as claimed.\\[12pt]
(ii) Fix $0<a<b$. We prove that there exists $C(a,b)>0$ such that
\[
A_N(b)\leq C(a,b)\bigl(N+A_N(a)\bigr)
\]
for every sufficiently large $N$. Dividing by $N$ and passing to the limit will then yield the result. First, we partition $[0,1)$ into 
\[
L_N=\left\lceil \frac{N}{a}\right\rceil
\]
intervals $J_1,\dots,J_{L_N}$ of equal length $\frac{1}{L_N}\leq \frac{a}{N}$. Let
\[
m_\ell=\#\{1\leq i\leq N:x_i\in J_\ell\}.
\]
so that $\sum_{\ell=1}^{L_N}m_\ell=N$. Since each interval $J_\ell$ has length at most $a/N$, any two points lying in the same $J_\ell$ are at distance at most $a/N$. As in (i) we thus obtain
\[
A_N(a) \geq \sum_{\ell=1}^{L_N}m_\ell(m_\ell-1).
\]
or equivalently 
\begin{align} \label{eq:ml}
\sum_{\ell=1}^{L_N}m_\ell^2 \leq A_N(a)+N.
\end{align}
If $x_i\in J_\ell$ and $x_j\in J_k$ with
\[
\norm{x_i-x_j}\leq \frac{b}{N},
\]
then for the cyclic (calculated modulo $L_N$) distance $d(\ell,k)$ between the interval indices $\ell$ and $k$ we have
\[
(d(\ell,k)-1)\frac{1}{L_N}\leq \frac{b}{N}.
\]
or in other words
\[
d(\ell,k)\leq \frac{bL_N}{N}+1.
\]
Since $\frac{L_N}{N}\to \frac{1}{a}$ as $N \to \infty$, it follows that, for all sufficiently large $N$, we have
\[
d(\ell,k)\leq \left\lfloor \frac{b}{a}\right\rfloor+2 =:K.
\]
Thus,
\[
A_N(b) \leq \sum_{\ell=1}^{L_N} \sum_{d(\ell,k)\leq K} m_\ell m_k \leq
\sum_{\ell=1}^{L_N} \sum_{d(\ell,k)\leq K} \frac{m_\ell^2+m_k^2}{2}.
\]
For each fixed $\ell$, there are at most $2K+1$ indices $k$ satisfying $d(\ell,k)\leq K$. Hence each $m_\ell^2$ appears with total multiplicity at most $2K+1$, and therefore
\[
A_N(b) \leq (2K+1)\sum_{\ell=1}^{L_N}m_\ell^2 \leq (2K+1)\left(N+A_N(a)\right).
\]
by \eqref{eq:ml}. Thus, with $C(a,b)$ as in the assertion of the theorem, we have
\[
A_N(b)\leq C(a,b)\bigl(N+A_N(a)\bigr)
\]
for all sufficiently large $N$. Dividing by $N$, we find and passing to the limit yield the claim
\[
f(b)\leq C(a,b)\bigl(1+f(a)\bigr).
\]
(iii) We choose $a>0$ with $f(a)=0$ and fix $b>0$. As $R_N(a) \to 0$, we see that $A_N(a)=o(N)$. Let $E_N(a)$ be the set of unordered pairs $i \neq j$ with $\norm{x_i-x_j}\leq \frac{a}{N}$. Then $\#E_N(a)=\frac{1}{2} A_N(a)=o(N)$. Choose a set $D_N\subset\{1,\dots,N\}$ containing at least one endpoint from each pair in $E_N(a)$ and let $
G_N=\{1,\dots,N\}\setminus D_N$. By construction, all points $x_i,x_j$ with $i, j\in G_N$ and $i \neq j$ satisfy
\[
\norm{x_i-x_j} > \frac{a}{N}.
\]
Now, we estimate the number of ordered pairs counted by $A_N(b)$ by splitting
them into three classes:
\[
G_N\text{--}G_N,\qquad D_N\text{--}G_N, \qquad D_N\text{--}D_N.
\]
We will argue that only the $G_N\text{--}G_N$ pairs contribute to $A_N(b)$ in a significant way. First, we consider the $G_N\text{--}G_N$ pairs. We may without loss of generality assume that the points $\{x_i:i\in G_N\}$ are ordered cyclically. Since $\norm{x_i-x_j} > a/N$ in this case, the $r$-th next point in either cyclic direction is at distance strictly larger than $ra/N$. Hence if $r\geq \lfloor \frac{b}{a} \rfloor +1$, then $\frac{ra}{N} > \frac{b}{N}$. Thus,  each point of $G_N$ has at most $\lfloor \frac{b}{a} \rfloor$ neighbors to the right and left respectively  within distance $b/N$. It follows, that the number of $G_N\text{--}G_N$ pairs counted by $A_N(b)$ is bounded from above by $2\lfloor \frac{b}{a}\rfloor N$.

Next, we consider pairs with one endpoint in $D_N$ and one endpoint in $G_N$. Since the points in $G_N$ have distance at least $a/N$, any interval of length $2b/N$ contains at most $C_1(a,b)=\left\lfloor \frac{2b}{a}\right\rfloor+2$ points. In other words, each point indexed by $D_N$ is within distance $b/N$ of at most $C_1(a,b)$ points indexed by $G_N$. Therefore the total number of ordered pairs of type $D_N\text{--}G_N$ or $G_N\text{--}D_N$ counted by $A_N(b)$ is at most $2C_1(a,b) \#D_N = o(N)$.

It remains to estimate the number of $D_N\text{--}D_N$ pairs. We again partition $[0,1)$ into
$$
L_N=\left\lceil \frac{N}{a}\right\rceil
$$
intervals $J_1,\dots,J_{L_N}$ of equal length $1/L_N$. As in (ii), the cyclic distance between indices $\ell$ and $k$ such that $x_i \in J_\ell$ and $x_j \in J_k$ and $\norm{x_i-x_j}\leq \frac{b}{N}$ is bounded by $K = \lfloor \frac{b}{a} \rfloor +2$ and hence the number of $D_N\text{--}D_N$ pairs which contribute to $A_N(b)$ is at most
\[
(2K+1)\bigl(A_N(a)+\#D_N\bigr)=o(N).
\]
Combining the three estimates for the different types of pairs, we obtain
\[
A_N(b) \leq 2  \left\lfloor \frac{b}{a} \right\rfloor N+o(N).
\]
Dividing by $N$ yields
$$
\limsup_{N\to\infty}R_N(b) \leq 2\left\lfloor \frac{b}{a}\right\rfloor .
$$
and hence
$$
f(b)\leq 2\left\lfloor \frac{b}{a}\right\rfloor.
$$
(iv) First we prove the sharpness of the lower bound. We build the sequence in blocks. To this end, let $\mathcal{B}_0 = \{ 0.5 \}$ be the first block, set $S_0=1$ and suppose inductively that the first $S_{m-1}$ points have already been chosen. We choose $M_m$ so large that
\[
\frac{S_{m-1}^2}{M_m} \to 0,
\]
e.g., $M_m = S_{m-1}^3$. Let 
\[
B_m=\lceil sM_m\rceil .
\]
be the size of the $m$-th block $\mathcal{B}_m$ and write
\[
B_m=q_m M_m+r_m,
\]
with $0\leq r_m<M_m$. We note that $q_m=\left\lfloor \frac{B_m}{M_m}\right\rfloor$ and $\lim_{m \to \infty} \frac{B_m}{M_m} = s$. The $m$-th block $\mathcal{B}_m$ will consist of $B_m$ points placed near the equidistant points
\[
\left\{0,\frac{1}{M_m},\dots,\frac{M_m-1}{M_m}\right\}.
\]
More precisely, we choose a shift $\tau_m\in\mathbb T$ avoiding coincidences with the previously chosen finitely many points. At $r_m$ of the equidistant points we place $q_m+1$ very close points, and at the remaining $M_m-r_m$ lattice sites we place $q_m$ very close points. The diameter of each cluster is chosen to be $\varepsilon_m>0$ so small that
\[
\varepsilon_m=o\left(\frac{1}{S_{m-1}+B_m}\right)
\]
and so large that the distance between all points in the distinct clusters in the $m$-th block is strictly greater than
\[
\frac{s}{S_{m-1}+B_m}.
\]
Indeed, this is possible because $\frac{1}{M_m}>\frac{s}{S_{m-1}+B_m}$ for all sufficiently large $m$. From now on, we restrict ourselves to the subsequence $N_m=S_{m-1}+B_m$, i.e., the length of the sequence after the $m$-th block has been appended. The aim is to show that $\lim_{m \to \infty} R_{N_m}(s)\longrightarrow G(s)$. Similarly as in the proof of (iii) there are three types of ordered pairs among the first $N_m$ points. First, old--old pairs, i.e., points with indices $\leq N_{m-1}$, contribute at most $\frac{S_{m-1}^2}{N_m}=o(1)$. 
Second we come to pairs with one old point and one new point. We claim that their contribution is $o(1)$ as well. Fix an old point $z$, i.e. one of the first $S_{m-1}$ points. We estimate how many points of the $m$-th block can lie in the interval
\[
\left\{y\in\mathbb [0,1): \norm{y-z}\leq \frac{s}{N_m}\right\}.
\]
The clusters of the $m$-th block are centered at the shifted points
\[
\tau_m+\frac{j}{M_m}, \qquad 0\leq j\leq M_m-1.
\]
and the spacing between neighboring cluster centers is $1/M_m$. Since every cluster has diameter at most $\varepsilon_m$, then the center itself lies within distance $
\frac{s}{N_m}+\varepsilon_m$ of $z$, if a point of the cluster centered at $\tau_m+j/M_m$ lies within distance $s/N_m$ of $z$. Hence the number of relevant cluster centers is bounded by
\[
2M_m\left(\frac{s}{N_m}+\varepsilon_m\right)+2,
\]
where the additional $+2$ accounts for boundary effects. Since $sM_m/N_m \leq 1$ and $
\varepsilon_m=o\left(\frac{1}{N_m}\right)$, it follows that
\[
M_m\varepsilon_m = o\left(\frac{M_m}{N_m}\right) = O_s(1).
\]
Consequently,
\[
2M_m\left(\frac{s}{N_m}+\varepsilon_m\right)+2 = O_s(1).
\]
Thus each old point can see only $O_s(1)$ clusters of the new block. Since the number of old elements divided by $N_m$ is of order $o(1)$, because $S_{m-1}^2/N_m \to 0$, it follows that the old--new and new--old pairs contribute $o(1)$ to $R_{N_m}(s)$.

It remains to count new--new pairs. By construction, no pairs from distinct clusters in the $m$-th block contribute to $R_N(s)$ but only pairs inside the same cluster do. A cluster containing $q_m$ points contributes $q_m(q_m-1)$ ordered pairs, and a cluster containing $q_m+1$ points contributes $(q_m+1)q_m$ ordered pairs. Hence the number of counted ordered pairs inside the $m$-th block $\mathcal{B}_m$ is
\[
(M_m-r_m)q_m(q_m-1)+r_m(q_m+1)q_m.
\]
Therefore,
\[
R_{N_m}(s) = \frac{(M_m-r_m)q_m(q_m-1)+r_m(q_m+1)q_m}{N_m} +o(1).
\]
Since $S_{m-1}^2/B_m\to 0$, we may replace $N_m$ by $B_m$ to obtain
\[
R_{N_m}(s) = \frac{(M_m-r_m)q_m(q_m-1)+r_m(q_m+1)q_m}{B_m} +o(1) =
\frac{M_m q_m(q_m-1)+2q_m r_m}{B_m}
+o(1).
\]
Taking into account
\[
\frac{B_m}{M_m}\to s, \qquad q_m\to h = \lfloor s \rfloor, \qquad \frac{r_m}{M_m}\to s-h.
\]
we see that
\[
\lim_{m \to \infty} R_{N_m}(s) = \frac{h(h-1)+2h(s-h)}{s} = \frac{h(2s-h-1)}{s} = G(s).
\]
which means that
\[
\liminf_{N\to\infty}R_N(s)=G(s).
\]
This completes the proof of the lower bound.

Finally, we prove the sharpness of the upper bound in (iii). Let $q\geq 2$ be an integer, and let $(x_n)_{n\geq 1}$ be the van der Corput sequence in base $q$ indexed so that $x_{n+1}=\phi_q(n)$,  where $\phi_q$ denotes radical inverse in base $q$, see e.g. \cite{KN74}. We are interested in the subsequence $N_m=q^m$. The first $N_m$ points are exactly the equidistant points
\[
\left\{0,\frac{1}{N_m},\frac{2}{N_m},\dots,\frac{N_m-1}{N_m}\right\}.
\]
Now fix an integer $k\geq 1$ and $a > 0$ with $\frac{k}{k+1}<a<1$ so that $\left\lfloor \frac{k}{a}\right\rfloor=k$. Since $a<1$, no two distinct equidistant points are within distance
$a/N_m$. Hence $R_{N_m}(a)=0$
and
\[
\liminf_{N\to\infty}R_{N}(a)=0.
\]
On the other hand, for fixed $k$ and $m$ sufficiently large, each lattice point has exactly $h$ neighbors to the right and $k$ neighbors to the left within distance $k/N_m$. Thus the number of ordered pairs $i\neq j$, satisfying
\[
\norm{x_i-x_j}\leq \frac{k}{N_m}
\]
is exactly $2kN_m$. Thus we have
$$
R_{N_m}(k) = \frac{1}{N_m}\cdot 2kN_m = 2k.
$$
and 
$$
\limsup_{N \to \infty}R_{N}(k) \geq 2k = 2 \left\lfloor \frac{k}{a} \right\rfloor.
$$
This completes the proof. 
\begin{remark} 
The restriction $a<1$ is essential in the proof, because we use the van der Corput sequence. As the definition of $R_N(s)$ uses $\leq$, each point of the equidistant points has two neighbors at distance exactly $1/N_m$ and so $R_{N_m}(1)=2$ but not $0$ (this would violate the desired condition $\liminf_{N \to \infty} R_N(a) =0$).    
\end{remark}
\begin{remark} Let $0<\alpha<1$, and define the weak $N$-point correlation by
\[
R_N^{(\alpha)}(s) = \frac{1}{N^{2-\alpha}} \#\left\{ 1\leq i\neq j\leq N \, : \,
\norm{x_i-x_j}\leq \frac{s}{N^\alpha} \right\},
\]
compare \cite{NP07}. The same arguments as in the proof Theorem~\ref{thm:main} yield analogous necessary bounds. More precisely, if $\lim_{N \to \infty} R_N^{(\alpha)}(s) = f_\alpha(s)$ for every $s\geq 0$, then $f_\alpha(s)\geq s$. Indeed, partitioning $[0,1)$ into $\lceil N^\alpha/s\rceil$ intervals of length at most $s/N^\alpha$ gives
$$
R_N^{(\alpha)}(s) \geq \frac{N^\alpha}{\lceil N^\alpha/s\rceil} - N^{\alpha-1},
$$
and the right-hand side tends to $s$. However, this bound is weaker than the one from \cite{Wei26d}. Moreover, for every $0<a<b$, the same interval argument as used in the proof of (ii) implies
\[
R_N^{(\alpha)}(b) \leq C(a,b)\left(N^{\alpha-1}+R_N^{(\alpha)}(a)\right)
\]
with 
\[
C(a,b)=2\left(\left\lfloor \frac{b}{a}\right\rfloor+2\right)+1.
\]
By passing to the limit we obtain $f_\alpha(b)\leq C(a,b)f_\alpha(a)$. Note that, since $f_\alpha(a)\geq a>0$, the analogue of the condition $f(a)=0$ in (iii) cannot occur for the weak $N$-point correlation.
\end{remark}

\paragraph{\textbf{AI declaration.}} In the course of preparing this paper, the author made use of ChatGPT 5.5 for idea generation, discussion of proof strategies, and language refinement. The author subsequently verified all claims, supplied missing details, and corrected inaccuracies. No AI-generated text or argument was included without critical review. The author assumes full responsibility for the accuracy and originality of the final manuscript.

\bibliographystyle{alpha}
\bibdata{references}
\bibliography{references}

\end{document}